\documentclass[11pt]{amsart}
\usepackage{amsmath,amsfonts,amssymb,amsthm,amstext,pgf,graphicx,hyperref,verbatim,lmodern,textcomp,color,young,tikz}
\usetikzlibrary{decorations}
\usepackage[mathscr]{euscript}
\usetikzlibrary{decorations.markings}
\usetikzlibrary{arrows}
\usetikzlibrary{decorations.pathreplacing,decorations.markings}
\theoremstyle{plain}
\newtheorem{theorem}{Theorem}[section]
\newtheorem{lemma}[theorem]{Lemma}

\theoremstyle{definition}

\newtheorem{exmp}{Example}[section]

\title{}
\begin{document}
\title [A graph-theoretic proof of Cramer's rule] {A graph-theoretic proof of Cramer's rule} 
\author[Sudip Bera]{Sudip Bera}
\address[Sudip Bera]{Faculty of Mathematics, Dhirubhai Ambani University, Gandhinagar-382421}
\email{sudip\_bera@dau.ac.in}
\keywords{Cramer's rule; combinatorial proof; determinants; digraphs}
\subjclass[2010]{05A19; 05A05; 05C30; 05C38}
\maketitle

\begin{abstract}
This note contains a new combinatorial proof of Cramer’s rule based on the Gessel-Viennot-Lindström Lemma.
\end{abstract}
\section{Introduction}
This paper presents a combinatorial proof of Cramer's rule. Such a proof offers a greater understanding of the underlying reasons for the validity of the result, rather than merely explaining the methodology \cite{14, 09,10}. Numerous concise proofs of Cramer’s rule are available on Wikipedia and its associated references \cite{Cramer-rule-six-proofs, Cramer-rule-Jacobi, Cramer-rule-Muir}. 

The rule was first published by Gabriel Cramer (1704--1752) in Appendix I of his \emph{Introduction à l’analyse des lignes courbes algébriques} \cite{Cramer-rule}, pages 657-659. While Theorem 1.1 is sometimes misattributed-Boyer, Hedman, and others suggest that Colin Maclaurin (1698-1746) was already aware of it by 1729 and included it in his posthumous \emph{Treatise of Algebra} (1748) \cite{Cramer-rule-Boyer,Cramer-rule-Hedman}. As a matter of fact, both Cramer and Maclaurin explicitly solved the $3\times3$ case, expressing each unknown as a ratio of two sums of six terms. They then sketched how these formulas extend to larger systems; neither, however, used the modern determinant concept, which emerged only in 1771 with Vandermonde \cite{Cramer-rule-Vandermonde}.

Furthermore, as observed in \cite{Cramer-rule-Kosinsky}, Maclaurin’s method for assigning signs to each summand is flawed. By contrast, Cramer’s approach-determining signs via the parity of the associated permutation is correct. Hence, the rule rightfully bears his name. In 1841, Carl Gustav Jacobi (1804-1851) introduced the first formal proof of Cramer's rule in his paper \cite{Cramer-rule-Jacobi}. However, this is not the earliest known demonstration; in 1825, Heinrich Ferdinand Scherk (1798-1885) published a 17-page inductive proof on the number of unknowns, outlined in \cite{Cramer-rule-Muir}. Recently, Doron Zeilberger provided a fully combinatorial proof in \cite{Cramer-rule-com-prf-zeilberger}. This paper presents a combinatorial proof of Cramer's rule utilizing the Gessel-Viennot-Lindström Lemma. 

Let $\Gamma$ represent a weighted, acyclic directed graph. Consider $P_1$ as a directed path from vertex $X$ to vertex $Y$ within $\Gamma$, and $P_2$ as another path extending from $Y$ to $Z$. The concatenation of the two paths, $P_1$ and $P_2$, is denoted as $P_1\bigodot P_2$, which traverses from vertex $X$ to vertex $Z$. A directed edge is represented by the initial vertex $U$ and the terminal vertex $V$ as $\overrightarrow{UV}.$ Let $A$ and $B$ be two fixed subsets of $V(\Gamma)$ both of cardinality $n$ respectively called set of \emph{initial vertices} and set of \emph{final vertices}, where $V(\Gamma)$ is the vertex set of the graph $\Gamma$. To these sets, we associate the \emph{path matrix} $M_{AB}=(m_{ij})_{n\times n}$, where $m_{ij}=\sum\limits_{P:A_i\rightarrow B_j}w(P)$, with $w(P)$ representing the product of the weights of all edges in the path $P.$ The notation $P:A_i\rightarrow B_j$ signifies a directed path that initiates at the vertex $A_i$ and concludes at the vertex $B_j.$ A \emph{path system} $\mathcal{P}$ from $A$ to $B$ consists of a permutation $\sigma$ and $n$ paths $P_i: A_i\rightarrow B_{\sigma(i)}$, with $\text{sgn} (\mathcal{P})=\text{sgn} (\sigma)$. The \emph{weight} of $\mathcal{P}$ is defined as $w(\mathcal{P})=\prod_{i=1}^nw(P_i)$. We refer to the path system as \emph{vertex-disjoint} if no two paths share a common vertex. Let $VD(\Gamma)$ denote the collection of vertex-disjoint path systems. It is straightforward to observe that ${\det}(M_{AB})=\sum\limits_{\mathcal{P}}\text{sgn}(\mathcal{P})w(\mathcal{P})$. However, the Gessel-Viennot-Lindström Lemma provides additional insights.

\begin{lemma}[Gessel-Viennot-Lindstr\"om \cite{15}]\label{Lemma:LGV-lemma} Let $\Gamma$ be a weighted, acyclic digraph and $M_{AB}$ be the path matrix of $\Gamma$. Then
${\det}(M_{AB})=\sum\limits_{\mathcal{P}\in VD(\Gamma)}\text{sgn}(\mathcal{P})w(\mathcal{P})$.
\end{lemma}
Note that the sum is $0$ if no path system exists from $A$ to $B$.
%Now we give an almost visual proof of Cramer's rule for solving the system of linear equations. Consider the system of linear equations:
%\begin{align*}
%a_{11}x_1+a_{12}x_2+\cdots +a_{1n}x_n=&b_1\\
%a_{21}x_1+a_{22}x_2+\cdots +a_{2n}x_n=&b_2 \\
%\vdots\hspace{1.5cm}\vdots\hspace{1cm}\vdots\hspace{1cm}\vdots\hspace{1cm}&\vdots\\
%a_{n1}x_1+a_{n2}x_2+\cdots +a_{nn}x_n=&b_n
%\end{align*}
%This system can be written as $AX=B$, where $A=(a_{ij})_{n\times n}$  is the $n$ ordered matrix,
%$X=(x_1, \cdots, x_n)^T$  and $B=(b_1, \cdots b_n)^T$ are two column vectors. Suppose $A_i (i=1, \cdots, n)$  is the matrix formed by replacing the $i$-th column of $A$ by the column vector $B$.
%\begin{theorem}[Cramer's rule \cite{Cramer-rule}]
%Let $AX=B$ be a system of  $n$ linear equations for $n$ unknown variables with $\det(A)\neq 0$. Then Cramer's rule says that  	
%\begin{equation*}
%x_i=\frac{\det(A_i)}{\det(A)}, (i=1, 2, \cdots, n)
%\end{equation*}
%\end{theorem}
We now present an almost visual demonstration of Cramer's rule for solving a system of linear equations. Consider the following system of equations:
\begin{align*}
a_{11}x_1+a_{12}x_2+\cdots +a_{1n}x_n=&b_1\\
a_{21}x_1+a_{22}x_2+\cdots +a_{2n}x_n=&b_2 \\
\vdots\hspace{.9cm}\vdots\hspace{1cm}\ddots\hspace{1cm}\vdots\hspace{1cm}&\vdots\\
a_{n1}x_1+a_{n2}x_2+\cdots +a_{nn}x_n=&b_n
\end{align*}
This system can be expressed in matrix form as $AX=B$, where $A=(a_{ij})_{n\times n}$ represents the $n \times n$ matrix, $X=(x_1, \cdots, x_n)^T$ is the column vector of the unknowns, and $B=(b_1, \cdots, b_n)^T$ is the column vector of constants. Let $A_i$ (for $i=1, \cdots, n$) denote the matrix obtained by substituting the $i$-th column of $A$ with the column vector $B$. 
\begin{theorem}[Cramer's rule \cite{Cramer-rule}]
For the system $AX=B$, consisting of $n$ linear equations with $n$ unknowns and $\det(A)\neq 0$, Cramer's rule states that
\begin{equation*}
x_i=\frac{\det(A_i)}{\det(A)}, \quad (i=1, \cdots, n).
\end{equation*}
\end{theorem}

\begin{figure}[ht!]
\tiny
%\tikzstyle{ver}=[]
%\tikzstyle{vert}=[circle, draw, fill=black!100, inner sep=0pt, minimum width=4pt]
%\tikzstyle{vertex}=[circle, draw, fill=black!00, inner sep=0pt, minimum width=4pt]
%\tikzstyle{edge} = [draw,thick,line width=.2 mm]
%\tikzstyle{edge_style} = [draw=black, line width=2, ultra thick]
%\tikzstyle{node_style} = [circle,draw=blue,fill=blue!20!,font=\sffamily\Large\bfseries]
\centering
\tikzset{->,>=stealth',auto, node distance=1cm,
	thick,main node/.style={circle,draw,font=\sffamily\Large\bfseries}}
\begin{tikzpicture}[scale=1.5]
%\tikzstyle{edge_style} = [draw=black, line width=2, ultra thick]
%\tikzstyle{node_style} = [circle,draw=blue,fill=blue!20!,font=\sffamily\Large\bfseries]
%\draw[->, line width=.4mm] (0,0) -- (1,0);	
\fill[blue!100!] (0, 2) circle (.05);
\fill[blue!100!] (0, -.05) circle (.05);
\node (A1) at (3,2)  {$\bf{\cdots}$};
\node (A1) at (3.3,2)  {$\bf{\cdots}$};
\node (A1) at (3.6,2)  {$\bf{\cdots}$};
%%%%%%
\node (A1) at (3,0)  {$\bf{\cdots}$};
\node (A1) at (3.3,0)  {$\bf{\cdots}$};
\node (A1) at (3.6,0)  {$\bf{\cdots}$};
%%%%%%%
\fill[blue!100!] (2, 2) circle (.05);
\fill[blue!100!] (2, -.05) circle (.05);
%%%%
\fill[blue!100!] (4.5, 0) circle (.05);
\fill[blue!100!] (7.5, 0) circle (.05);
%%%%%%%
\fill[blue!100!] (6, 2) circle (.05);
\fill[blue!100!] (6, 0) circle (.05);
%%%%%%%
\fill[blue!100!] (9, 2) circle (.05);
\fill[blue!100!] (9, 0) circle (.05);
\fill[blue!100!] (6, -2) circle (.06);
%%%%%%%%%%%%%%%
\node (A1) at (0,2.3)  {$\bf{A_1}$};
\node (A1) at (0,-.3)  {$\bf{B_1}$};
\node (A1) at (-.2,1)  {$\bf{a_{11}}$};
\node (A1) at (.6,1.8)  {$\bf{a_{12}}$};
\node (A1) at (1.5,1.8)  {$\bf{a_{21}}$};
\node (A1) at (2.3,1)  {$\bf{a_{22}}$};
\node (A1) at (4.9,1)  {$\bf{a_{ik}}$};
\node (A1) at (3.7,1)  {$\bf{a_{2k}}$};
\node (A1) at (5.8,1)  {$\bf{a_{ii}}$};
\node (A1) at (7,1)  {$\bf{a_{i\ell}}$};
\node (A1) at (9.25,1)  {$\bf{a_{nn}}$};
%%%%%%%%%%%%%%%%%
\node (A1) at (2,2.3)  {$\bf{A_2}$};
\node (A1) at (2,-.3)  {$\bf{B_2}$};
%%%%
\node (A1) at (6,2.3)  {$\bf{A_i}$};
\node (A1) at (5.8,-.3)  {$\bf{B_i}$};
\node (A1) at (7,2)  {$\bf{\cdots}$};
\node (A1) at (7.3,2)  {$\bf{\cdots}$};
\node (A1) at (7.6,2)  {$\bf{\cdots}$};
%%%%%%%
\node (A1) at (5.3, 0)  {$\bf{\cdots}$};
\node (A1) at (6.6,0)  {$\bf{\cdots}$};
\node (A1) at (8,0)  {$\bf{\cdots}$};
%%%%
\node (A1) at (4.5,-.3)  {$\bf{B_k}$};
\node (A1) at (7.5,-.3)  {$\bf{B_{\ell}}$};
%%%%
\node (A1) at (9,2.3)  {$\bf{A_n}$};
\node (A1) at (9,-.3)  {$\bf{B_{n}}$};
\node (A1) at (6,-2.3)  {$\bf{X}$};
%%%%
\draw[->,line width=.2 mm,red] (0,2) -- (0,0);
\draw[->, line width=.2 mm,red] (2,2) -- (2,0);
\draw[->, line width=.2 mm,red] (2,2) -- (4.5,0);
%%%%%%%%%%
\draw[->, line width=.2 mm,red] (0,2) -- (2,0);
\draw[->, line width=.2 mm,red] (2,2) -- (0,0);
%%%%%%%
\draw[->, line width=.2 mm,red] (6,2) -- (6,0);
%\draw[->, line width=.3 mm] (2,2) -- (0,0);
%%%
\draw[->, line width=.2 mm,red] (6,2) -- (4.5,0);
\draw[->, line width=.2 mm,red] (6,2) -- (7.5,0);
%%%%%
\draw[->, line width=.2 mm,red] (9,2) -- (9,0);
\draw[->, line width=.2 mm,red] (6,0) -- (6,-1.95);
\draw[->, line width=.2 mm,red] (0,0) -- (5.95,-2);
\draw[->, line width=.2 mm,red] (2,0) -- (5.95,-2);
\draw[->, line width=.2 mm,red] (4.5,0) -- (6,-1.95);
\draw[->, line width=.2 mm,red] (7.5,0) -- (6.05,-2);
\draw[->, line width=.2 mm,red] (9,0) -- (6.05,-2);
%%%%%%%%%%%%%%%%%%%%%%%%%%
\node (1) at (.7,-.1)  {$\bf{x_1}$};
\node (1) at (2.4,-.08)  {$\bf{x_2}$};
\node (1) at (4.9,-.3)  {$\bf{x_k}$};
\node (1) at (6.2,-.3)  {$\bf{x_i}$};
\node (1) at (7.1,-.3)  {$\bf{x_{\ell}}$};
\node (1) at (8.3,-.3)  {$\bf{x_{n}}$};

\end{tikzpicture}
\caption{
	$\Gamma$ is a weighted digraph having a directed edge from $A_i$ to $B_j$ with weight $a_{ij}$ for each $i, j \in [n]$.}

\label{fig:cramer-general}	
\end{figure}

\begin{proof}
Our objective is to demonstrate that $x_i \det(A) = \det(A_i)$ for every $i \in [n]$. Consider the directed graph $\Gamma$ illustrated in Figure \ref{fig:cramer-general}. The graph $\Gamma$ is a weighted digraph having directed edge from $A_i$ to $B_j$ with weight $a_{ij}$ for each $i, j\in [n]$ and the weight of the edge $\overrightarrow{B_iX}$ is $x_i,$ for each $i\in [n].$ Let $A = \{ A_1, \cdots, A_n \}$ represent the initial set of vertices, while $B = \{ B_1, \cdots, B_{i-1}, X, B_{i+1}, \cdots, B_n \}$ denotes the terminal set of vertices in $\Gamma.$ The weight associated with the edge connecting vertex $A_i$ to vertex $B_j$ in the graph $\Gamma$ is denoted as $a_{ij}.$ Furthermore, the weight of the edge from vertex $B_i$ to vertex $X$ is represented by $x_i$. It is important to note that
\[
\sum\limits_{P: A_j\rightarrow X} w(P) = \sum\limits_{k=1}^n a_{jk} x_k, \text{ for all } j \in [n].
\]
Consequently, the $i$-th column of the path matrix $M_{AB}$ in the graph $\Gamma$ can be expressed as follows:
\begin{align*}
\left(
\begin{array}{c}
\sum\limits_{k=1}^n{a_{1k}x_k} \\
\sum\limits_{k=1}^n{a_{2k}x_k} \\
\vdots \\
\sum\limits_{k=1}^n{a_{nk}x_k}
\end{array}
\right)=\left(\begin{array}{c}
b_1 \\
b_2\\
\vdots\\
b_n
\end{array}\right).
\end{align*}
Furthermore, it is evident that the column $C_j, j\in [n]\setminus \{i\}$ of the path matrix $M_{AB}$ is represented as:
\begin{align*}
\left(\begin{array}{c}
a_{1j} \\
a_{2j}\\
\vdots\\
a_{nj}
\end{array}\right).
\end{align*}
Thus, the path matrix $M_{AB}$ can be formulated as:
\begin{align*}
\left(
\begin{array}{ccccc}
a_{11}&\cdots& \sum\limits_{k=1}^n{a_{1k}x_k}&\cdots& a_{1n} \\
a_{21}&\cdots&\sum\limits_{k=1}^n{a_{2k}x_k}& \cdots& a_{2n}\\
\vdots&\ddots& \vdots&\ddots&\vdots\\
a_{(n-1)1}&\cdots&\sum\limits_{k=1}^n{a_{(n-1)k}x_k}&\cdots&a_{(n-1)n}\\
a_{n1}&\cdots&\sum\limits_{k=1}^n{a_{nk}x_k}&\cdots&a_{nn}
\end{array}
\right)=A_i.
\end{align*}
According to Lemma \ref{Lemma:LGV-lemma}, it follows that $\det(A_i)=\sum\limits_{\mathcal{P}\in VD(\Gamma)}\text{sgn}(\mathcal{P})w(\mathcal{P})$. From Figure \ref{fig:cramer-general}, it is evident that the set $\mathcal{P}=\{P_1, \cdots, P_n\}$ constitutes a vertex disjoint path system in the induced graph $\Gamma\setminus \{X\}$, with the initial vertex set being $\{A_1, \cdots, A_n\}$ and the terminal vertex set being $\{B_1, \cdots, B_n\}$ if and only if $\mathcal{\bar{P}}=\{P_1, \cdots, P_{i-1}, P_i\bigodot \overrightarrow{B_i X}, P_{i+1},\cdots, P_n\}$ forms a vertex disjoint path system in the graph $\Gamma$, where $A=\{A_1, \cdots, A_n\}$ and $B=\{B_1, \cdots, B_{i-1}, X, B_{i+1}, \cdots, B_n\}$ represent the initial and terminal vertex sets of $\Gamma$, respectively. Furthermore, it is important to observe that $w(\mathcal{\bar{P}})=x_iw(\mathcal{P})$ and $\text{sgn}(\mathcal{\bar{P}})=\text{sgn}(\mathcal{P})$. Consequently, we have
\begin{align*}
\left(\sum\limits_{\mathcal{P}\in VD(\Gamma)}\text{sgn} (\mathcal{P})w(\mathcal{P})\right)&=x_i \left(\sum\limits_{\mathcal{P}\in VD(\Gamma\setminus \{X\})}\text{sgn} (\mathcal{P})w(\mathcal{P})\right) \\
\Rightarrow \det(A_i)&=x_i\det(A).
\end{align*}
This concludes the proof.

\begin{exmp}
Here we explain the idea of the proof for the case $n=3.$ Consider the graph $\Gamma$ in Figure \ref{fig:cramer}.
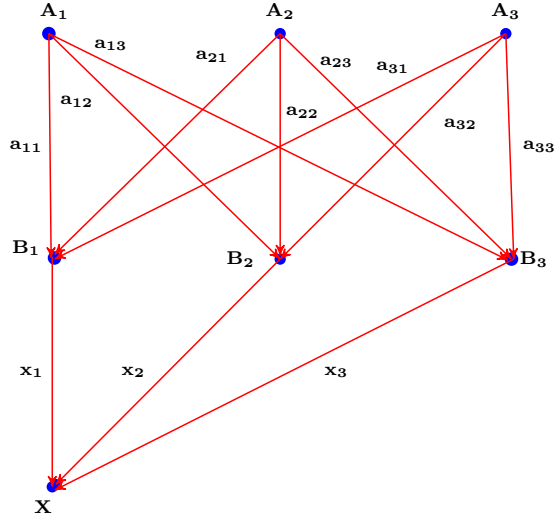
\begin{figure}[ht!]
\tiny
%\tikzstyle{ver}=[]
%\tikzstyle{vert}=[circle, draw, fill=black!100, inner sep=0pt, minimum width=4pt]
%\tikzstyle{vertex}=[circle, draw, fill=black!00, inner sep=0pt, minimum width=4pt]
%\tikzstyle{edge} = [draw,thick,line width=.2 mm]
%\tikzstyle{edge_style} = [draw=black, line width=2, ultra thick]
%\tikzstyle{node_style} = [circle,draw=blue,fill=blue!20!,font=\sffamily\Large\bfseries]
\centering
\tikzset{->,>=stealth',auto, node distance=1cm,
	thick,main node/.style={circle,draw,font=\sffamily\Large\bfseries}}
%\tikzset{->-/.style={decoration={
				%markings,
			%mark=at position #1 with {\arrow{>}}},postaction={decorate}}}

\begin{tikzpicture}[scale=1.5]
%\tikzstyle{edge_style} = [draw=black, line width=2, ultra thick]
%\tikzstyle{node_style} = [circle,draw=blue,fill=blue!20!,font=\sffamily\Large\bfseries]
%\draw[->, line width=.4mm] (0,0) -- (1,0);	

\fill[blue!100!] (2, 0) circle (.05);
\fill[blue!100!] (2, 2) circle (.05);
\fill[blue!100!] (-0.05,2) circle (.06);
\fill[blue!100!] (0.01,-2) circle (.05);
\fill[blue!100!] (-.02,-2.02) circle (.05);
\fill[blue!100!] (4.05,0) circle (.06);
\fill[blue!100!] (4,2) circle (.05);
\fill[blue!100!] (0,0.01) circle (.06);
%%%%%%%%%%%%%%%%%%%%%%%%%%%%%%%%%%%%%
\node (A1) at (0,2.2)  {$\bf{A_1}$};
\node (A3) at (4,2.2)  {$\bf{A_3}$};
\node (B1) at (-.25,.1)   {$\bf{B_1}$};
\node (B2) at (1.65, 0)   {$\bf{B_2}$};
\node (B3) at (4.25, 0)   {$\bf{B_3}$};
\node (X) at (-.1,-2.2)   {$\bf {X}$};
\node (A2) at (2,2.2)  {$\bf{A_2}$};
%\draw[->, line width 1mm] (4,2)--(2,0);

\draw[->, line width=.2 mm,red] (4,2) -- (2,0);
\draw[->, line width=.2 mm,red] (2,0)--(0.01,-2);
\draw[->, line width=.2 mm,red] (-0.05,2)--(-.03,0.01);
\draw[->, line width=.2 mm,red] (-.02,0.02)--(-.02,-2.02);
\draw[->, line width=.2 mm,red] (2,2)--(4.05,0);
\draw[->, line width=.2 mm,red] (-0.05,2)-- (4.02,-.01);
\draw[->,line width=.2 mm,red] (4,2)--(4.07,0);
\draw[->, line width=.2 mm,red] (4,2)--(0,0);
\draw[->, line width=.2 mm,red] (2,2)--(2,0.03);
\draw[->,line width=.2 mm,red] (4.07,0)-- (0,-2.05);
\draw[->, line width=.2 mm,red] (-0.05,2)--(2,0);
\draw[->, line width=.2 mm,red] (2,2)--(0,0.04);

%\draw [fill=orange] (0.1,0.1) rectangle (0.2,0.2);

%\draw[->-=.8] (4,2) -- (2,0);
%\draw[->-=.9] (4,2)--(0,-2);
%\draw[->-=.8] (0,2)--(0,0);
%\draw[->-=.8] (0,0)--(0,-2);
%\draw[->-=.8] (2,2)--(4,0);
%\draw[->-=.9] (0,2)-- (4,0);
%\draw[->-=.8] (4,2)--(4,0);
%\draw[->-=.8] (4,2)--(0,0);
%\draw[->-=.8] (2,2)--(2,0);
%\draw[->-=.9] (4,0)-- (0,-2);
%\draw[->-=.8] (0,2)--(2,0);
%\draw[->-=.8] (2,2)--(0,0);

\node (1) at (-.2,-1)  {$\bf{x_1}$};
\node (1) at (.7,-1)  {$\bf{x_2}$};
\node (1) at (2.5,-1)  {$\bf{x_3}$};
\node (1) at (-.25,1)  {$\bf{a_{11}}$};
\node (1) at (.2,1.4)  {$\bf{a_{12}}$};
\node (1) at (.5,1.9)  {$\bf{a_{13}}$};
\node (1) at (1.4,1.8) {$\bf{a_{21}}$};
\node (1) at (2.2,1.3) {$\bf{a_{22}}$};
\node (1) at (2.5,1.75)  {$\bf{a_{23}}$};
\node (1) at (3,1.7)  {$\bf{a_{31}}$};
\node (1) at (3.6,1.2)  {$\bf{a_{32}}$};
\node (1) at (4.3,1)  {$\bf{a_{33}}$};
\end{tikzpicture}
\caption{$\Gamma$ is a weighted digraph having edge weight $a_{ij}$ for each directed  edge $A_i$ to $B_j$ and $x_i$ for each edge $B_i$ to $X.$}
\label{fig:cramer}	
\end{figure}    
\end{exmp}
We aim to demonstrate that $\det(A_1)=x_1\det(A)$. Let us define the sets $A=\{ A_1, A_2, A_3\}$ and $B=\{ X, B_2, B_3\}$ as the initial and terminal sets of vertices in the graph $\Gamma$, respectively. It is straightforward to observe that $w(\mathcal{\bar{P}})=x_1w(\mathcal{P})$ and $\text{sgn}(\mathcal{\bar{P}})=\text{sgn}(\mathcal{P})$, where $\mathcal{\bar{P}}$ and $\mathcal{P}$ represent vertex-disjoint path systems in the graphs $\Gamma$ and $\Gamma\setminus\{X\}$, respectively. Consequently, we have the following relationship:
\begin{align*}
\left(\sum\limits_{\mathcal{P}\in VD(\Gamma)}\text{sgn} (\mathcal{P})w(\mathcal{P})\right)&=x_1 \left(\sum\limits_{\mathcal{P}\in VD(\Gamma\setminus \{X\})}\text{sgn} (\mathcal{P})w(\mathcal{P})\right) \\
\Rightarrow \det(A_1)&=x_1\det(A).
\end{align*}
\end{proof}

\subsection*{Data availability statement} Availability of data and materials are not applicable. 
\subsection*{Conflict of Interest}The author does not have disclosed any competing interests.	
\subsection*{Acknowledgement} The authors wish to sincerely thank the referee for her/his comments and suggestions, thus improving the submitted version of the paper.

\bibliographystyle{amsplain}
\bibliography{gen-inv-lcp5}

\end{document}